\documentclass[11pt,a4paper,fleqn]{article}
\usepackage{a4wide,amsfonts,amsmath,latexsym,amssymb,euscript,graphicx,units,mathrsfs}

\usepackage{graphicx}
\usepackage{color}
\usepackage{amssymb}
\usepackage{amssymb}
\usepackage[T1]{fontenc}
\usepackage{latexsym}
\usepackage{xypic}
\usepackage{eufrak}
\usepackage{euscript}
\usepackage{amsfonts,amsmath}
\usepackage{verbatim}
\usepackage{fancyhdr}
\usepackage[english]{babel}
\usepackage{mathrsfs}
\usepackage{units}

\newtheorem{prop}{Proposition}[section]
\newtheorem{proposition}[prop]{Proposition}
\newtheorem{cor}[prop]{Corollary}
\newtheorem{lemma}[prop]{Lemma}
\newtheorem{rem}[prop]{Remark}
\newtheorem{remark}[prop]{Remark}

\newtheorem{theorem}[prop]{Theorem}
\newtheorem{thm}[prop]{Theorem}

\renewcommand{\geq}{\geqslant}
\def\leq{\leqslant}

\newcommand{\R}{\mathbb{R}}

\def\HH{\EuFrak H}

\def\1{{\mathbf{1}}}

\def\1{{\mathbf{1}}}
\def\0.5{{\frac{1}{2}}}

\newcommand{\qed}{\nopagebreak\hspace*{\fill}
{\vrule width6pt height6ptdepth0pt}\par}

\begin{document}

\begin{center}
{\Large{\bf Second order Poincar\'e inequalities and CLTs on Wiener space}}
\normalsize
\\~\\ by Ivan Nourdin\footnote{Laboratoire de Probabilit\'es et Mod\`eles Al\'eatoires, Universit\'e
 Pierre et Marie Curie (Paris VI), Bo\^ite courrier 188, 4 place Jussieu, 75252 Paris Cedex 05, France. Email: \texttt{ivan.nourdin@upmc.fr}},
 Giovanni Peccati\footnote{Equipe Modal'X, Universit\'{e} Paris Ouest -- Nanterre la D\'{e}fense, 200 Avenue de la République, 92000 Nanterre, and LSTA, Universit\'{e} Paris VI, France. Email: \texttt{giovanni.peccati@gmail.com}} and Gesine Reinert\footnote{Department of Statistics, University of Oxford, 1 South Parks Road, Oxford OX1 3TG, UK. Email \texttt{reinert@stats.ox.ac.uk}}\\ {\it Universit\'e Paris VI, Universit\'e Paris Ouest and Oxford University}\\~\\
\end{center}

{\small \noindent {\bf Abstract:} We prove infinite-dimensional second order Poincar\'{e} inequalities on Wiener space, thus closing a circle of ideas linking limit theorems for functionals of Gaussian fields, Stein's method and Malliavin calculus. We provide two applications: (i) to a new ``second order'' characterization of CLTs on a fixed Wiener chaos, and (ii) to linear functionals of Gaussian-subordinated fields.}
\\



\section{Introduction }
Let $N\sim\mathscr{N}(0,1)$ be a standard Gaussian random variable. In its most basic formulation, the
 {\sl Gaussian Poincaré inequality} states that, for every differentiable function $f:\R\to\R$,
\begin{equation}\label{PoincBasic}
{\bf Var}f(N) \leq Ef'(N)^2,
\end{equation}
with equality if and only if $f$ is affine. The estimate (\ref{PoincBasic}) is a fundamental tool of stochastic analysis:
it implies that, if the random variable $f'(N)$ has a small $L^2(\Omega)$ norm, then $f(N)$ has necessarily small
fluctuations. Relation (\ref{PoincBasic}) has been first proved by Nash in \cite{Nash}, and then rediscovered by
Chernoff in \cite{Chernoff} (both proofs use Hermite polynomials). The Gaussian Poincaré inequality admits extensions
in several directions, encompassing both the case of smooth functionals of multi-dimensional (and possibly infinite-dimensional)
 Gaussian fields, and of non-Gaussian probability distributions -- see e.g. Bakry {\it et al.} \cite{BBCG},
Bobkov \cite{Bobkov}, Cacoullos {\it et al.}, Chen \cite{Chen1985, Chen1987, Chen1988}, Houdré and Perez-Abreu \cite{HouPA},
and the references therein. In particular, the results proved in \cite{HouPA} (which make use of the Malliavin calculus)
allow to recover the following infinite-dimensional version of (\ref{PoincBasic}). Let $X$ be an isonormal Gaussian process
over some real separable Hilbert space $\mathfrak{H}$ (see Section \ref{S : Pre}), and let $F\in\mathbb{D}^{1,2}$ be a
Malliavin-differentiable functional of $X$. Then, the Malliavin derivative of $F$, denoted by $DF$, is a random element
with values in $\mathfrak{H}$, and it holds that
\begin{equation}\label{INFpoinc}
{\bf Var} F \leq E\|DF\|^2_{\mathfrak{H}},
\end{equation}
with equality if and only if $F$ has the form of a constant plus an element of the first Wiener chaos of $X$. In Proposition
\ref{nicelemma} below we shall prove a more general version of (\ref{INFpoinc}), involving central moments of arbitrary even
orders and based on the techniques developed in \cite{stein-ptrf}. Note that (\ref{INFpoinc}) contains as a special case the
well-known fact that, if $F = f(X_1,...,X_d)$ is a smooth function of  i.i.d. $\mathscr{N}(0,1)$ random variables $X_1,...,X_d$,
then
\begin{equation}
\label{vectorPOINC}
{\bf Var} F \leq E \|\nabla f(X_1,...,X_d) \|^2_{\mathbb{R}^d},
\end{equation}
where $\nabla f$ is the gradient of $f$.

Now suppose that the random variable $F =f(X_1,...,X_d)$ (where the $X_1,...,X_d$ are again i.i.d. $\mathscr{N}(0,1)$) is
such that $f$ is twice differentiable. In the recent paper \cite{Chatterjee_ptrf}, Chatterjee has pointed out that if one
focuses also on the $d\times d$ Hessian matrix ${\rm Hess}\, f$, and not only on $\nabla f$, then one can state an inequality
assessing the {\sl total variation distance} (see Section \ref{SS :distances}, \eqref{TV dst}) between the law of $F$ and the law of a
Gaussian random variable with matching mean and variance. The precise result goes as follows
(see \cite[Theorem 2.2]{Chatterjee_ptrf}). Let $E(F) = \mu$, ${\bf Var}F = \sigma^2>0$,
$Z\sim \mathscr{N}(\mu,\sigma)$, and denote by $d_{TV}(F,Z)$ the total variation distance between the laws
of $F$ and $Z$, see (\ref{TV dst}). Then
\begin{equation}\label{2ndChat}
d_{TV}(F,Z) \leq \frac{2\sqrt{5}}{\sigma^2} E[\|{\rm Hess}\, f(X_1,...,X_d)\|_{op}^4]^{\frac14} \times E[\|\nabla f(X_1,...,X_d)\|_{\mathbb{R}^d}^4]^{\frac14},
\end{equation}
where $\|{\rm Hess}\, f(X_1,...,X_d)\|_{op}$ is the operator norm of the (random) matrix ${\rm Hess}f(X_1,...,X_d)$. A relation
such as (\ref{2ndChat}) is called a {\sl second order Poincaré inequality}: it is proved in \cite{Chatterjee_ptrf} by combining
(\ref{vectorPOINC}) with an adequate version of {\sl Stein's method} (see e.g. \cite{Chen_Shao_sur, Reinert_sur}).

In \cite[Remark 3.6]{stein-ptrf} the first two authors of the present paper pointed out that the finite-dimensional Stein-type
inequalities leading to Relation (\ref{2ndChat}) are special instances of much more general estimates, which  can be obtained by
 combining Stein's method and Malliavin calculus on an infinite-dimensional Gaussian space. It is therefore natural to ask
whether the results of \cite{stein-ptrf} can be used in order to obtain a general version of (\ref{2ndChat}), involving a
``distance to Gaussian'' for smooth functionals of arbitrary infinite-dimensional Gaussian fields. We shall show that the
answer is positive. Indeed, one of the principal achievements of this paper is the proof of the following
statement ($d_W$
denotes the Wasserstein distance, see \eqref{W dst}):

\begin{theorem}\label{thm1}
{\bf (Second order infinite-dimensional Poincar\'{e} inequality)}
Let $X$ be an isonormal Gaussian process over some real separable Hilbert space $\mathfrak{H}$, and let $F\in \mathbb{D}^{2,4}$.
Assume that  $E\left(
F\right) =\mu $ and ${\rm Var}\left( F\right) =\sigma ^{2}>0$.
Let $Z\sim \mathscr{N}( \mu ,\sigma ^{2})$. Then
\begin{equation}
d_{W}\left( F,Z \right) \leq \frac{\sqrt{10}}{2 \sigma} E\left[ \Vert D^{2}F\Vert _{op}^{4}\right]
^{\frac14}\times E\left[ \left\Vert DF\right\Vert _{
\mathfrak{H}
}^{4}\right] ^{\frac14}.  \label{Poincw}
\end{equation}
If, in addition,
the law of $F$ is absolutely continuous with respect to the Lebesgue measure, then
\begin{equation}
d_{TV}\left( F,Z \right) \leq \frac{\sqrt{10}}{\sigma ^{2}} E\left[ \Vert D^{2}F\Vert _{op}^{4}\right]
^{\frac14}\times E\left[ \left\Vert DF\right\Vert _{
\mathfrak{H}
}^{4}\right] ^{\frac14}.   \label{Poinc}
\end{equation}
\end{theorem}

The class $\mathbb{D}^{2,4}$ of twice Malliavin-differentiable functionals is formally defined in Section \ref{S : Pre}; note
 that $D^2F$ is a random element with values in $\mathfrak{H}^{\odot 2}$ (the symmetric tensor product of $\mathfrak{H}$ with
itself) and that we used $\left\Vert D^{2}F\right\Vert _{op}$ to indicate the operator norm (or, equivalently, the
spectral radius) of the random Hilbert-Schmidt operator $f\mapsto \langle f, D^2 F \rangle_\mathfrak{H}$. The proof of
Theorem \ref{thm1} is detailed in Section \ref{SS :ProofT1}.
As discussed in Section \ref{SS : CONTR}, a crucial point is that Theorem \ref{thm1} leads to further (and very useful)
inequalities, which  we name {\sl random contraction inequalities}. These estimates involve a ``contracted version'' of the
second derivative $D^2F$, and will lead (see Section \ref{S : CharCLTWC}) to the proof of new necessary and sufficient
conditions which ensure that a sequence of random variables belonging to fixed Wiener chaos converges in law to a standard
Gaussian random variable. This result generalizes and unifies the findings contained in \cite{stein-ptrf, NO, NuPec, PT},
and virtually closes a very fruitful circle of recent ideas linking Malliavin calculus, Stein's method and central limit
theorems (CLTs) on Wiener space (see also \cite{NouPecNcltAOP}). The role of contraction inequalities is further explored
in Section \ref{S : GaussSub}, where we study CLTs for linear functionals of Gaussian subordinated fields.

\smallskip

The rest of the paper is organized as follows. In Section \ref{S : Pre} we recall some preliminary results involving Malliavin
operators. Section \ref{SS : TechLemmas} concerns Poincaré type inequalities and bounds on distances between probabilities.
Section \ref{S : MAIN} deals with the proof of Theorem \ref{thm1}, as well as with ``random contraction inequalities''. Section
\ref{S : CharCLTWC} and Section \ref{S : GaussSub} focus, respectively, on CLTs on Wiener chaos and on CLTs for Gaussian
subordinated fields. Finally, Section \ref{S : Generaliz} is devoted to a version of (\ref{Poincw}) for random
variables of the type $F= (F_1, \ldots, F_d)$.

\section{Preliminaries}\label{S : Pre}
We shall now present the basic elements of Gaussian analysis and Malliavin calculus that are used in this paper. The reader
is referred to the two monographs by Malliavin \cite{MallBook} and Nualart \cite{Nbook} for any unexplained definition or
 result.

\smallskip

Let $\EuFrak H$ be a real separable Hilbert space. For any $q\geq 1$ let $\EuFrak H^{\otimes q}$ be the $q$th tensor product of $\EuFrak H$ and denote
by $\EuFrak H^{\odot q}$ the associated $q$th symmetric tensor product. We write $X=\{X(h),h\in \EuFrak H\}$ to indicate
an isonormal Gaussian process over
$\EuFrak H$, defined on some probability space $(\Omega ,\mathcal{F},P)$.
This means that $X$ is a centered Gaussian family, whose covariance is given in terms of the
inner product of $\EuFrak H$ by $E\left[ X(h)X(g)\right] =\langle h,g\rangle _{\EuFrak H}$. We also assume that
 $\mathcal{F}$ is
generated by $X$.

For every $q\geq 1$, let $\mathcal{H}_{q}$ be the $q$th Wiener chaos of $X$,
that is, the closed linear subspace of $L^2(\Omega ,\mathcal{F},P)$
generated by the random variables of the type $\{H_{q}(X(h)),h\in \EuFrak H,\left\|
h\right\| _{\EuFrak H}=1\}$, where $H_{q}$ is the $q$th Hermite polynomial
defined as $H_q(x) = (-1)^q e^{\frac{x^2}{2}}
 \frac{d^q}{dx^q} \big( e^{-\frac{x^2}{2}} \big)$.
We write by convention $\mathcal{H}_{0} = \mathbb{R}$. For
any $q\geq 1$, the mapping $I_{q}(h^{\otimes q})=q!H_{q}(X(h))$ can be extended to a
linear isometry between \ the symmetric tensor product $\EuFrak H^{\odot q}$
equipped with the \ modified norm $\sqrt{q!}\left\| \cdot \right\| _{\EuFrak H^{\otimes q}}$ and the $q$th Wiener chaos $\mathcal{H}_{q}$. For $q=0$ we write $I_{0}(c)=c$, $c\in\mathbb{R}$.

It is well-known (Wiener chaos expansion) that $L^2(\Omega ,\mathcal{F},P)$
can be decomposed into the infinite orthogonal sum of the spaces $\mathcal{H}_{q}$. Therefore, any square integrable random variable
$F\in L^2(\Omega ,\mathcal{F},P)$ admits the following chaotic expansion
\begin{equation}
F=\sum_{q=0}^{\infty }I_{q}(f_{q}),  \label{E}
\end{equation}
where $f_{0}=E[F]$, and the $f_{q}\in \EuFrak H^{\odot q}$, $q\geq 1$, are
uniquely determined by $F$. For every $q\geq 0$ we denote by $J_{q}$ the
orthogonal projection operator on the $q$th Wiener chaos. In particular, if
$F\in L^2(\Omega ,\mathcal{F},P)$ is as in (\ref{E}), then
$J_{q}F=I_{q}(f_{q})$ for every $q\geq 0$.

Let $\{e_{k},\,k\geq 1\}$ be a complete orthonormal system in $\EuFrak H$.
Given $f\in \EuFrak H^{\odot p}$ and $g\in \EuFrak H^{\odot q}$, for every
$r=0,\ldots ,p\wedge q$, the \textit{contraction} of $f$ and $g$ of order $r$
is the element of $\EuFrak H^{\otimes (p+q-2r)}$ defined by
\begin{equation}
f\otimes _{r}g=\sum_{i_{1},\ldots ,i_{r}=1}^{\infty }\langle
f,e_{i_{1}}\otimes \ldots \otimes e_{i_{r}}\rangle _{\EuFrak H^{\otimes
r}}\otimes \langle g,e_{i_{1}}\otimes \ldots \otimes e_{i_{r}}
\rangle_{\EuFrak H^{\otimes r}}.  \label{v2}
\end{equation}
Notice that $f\otimes _{r}g$ is not necessarily symmetric: we denote its
symmetrization by $f\widetilde{\otimes }_{r}g\in \EuFrak H^{\odot (p+q-2r)}$.
Moreover, $f\otimes _{0}g=f\otimes g$ equals the tensor product of $f$ and
$g$ while, for $p=q$, $f\otimes _{q}g=\langle f,g\rangle _{\EuFrak H^{\otimes q}}$.
In the particular case where $\EuFrak H=L^2(A,\mathcal{A},\mu )$, where
$(A,\mathcal{A})$ is a measurable space and $\mu $ is a $\sigma $-finite and
non-atomic measure, one has that $\EuFrak H^{\odot q}=L_{s}^{2}(A^{q},
\mathcal{A}^{\otimes q},\mu ^{\otimes q})$ is the space of symmetric and
square integrable functions on $A^{q}$. Moreover, for every $f\in
\EuFrak H^{\odot q}$, $I_{q}(f)$ coincides with the multiple Wiener-It\^{o} integral
of order $q$ of $f$ with respect to $X$ introduced by It\^{o} in \cite{It}.
In this case, (\ref{v2}) can be written as
\begin{eqnarray*}
(f\otimes _{r}g)(t_1,\ldots,t_{p+q-2r})
&=&\int_{A^{r}}f(t_{1},\ldots ,t_{p-r},s_{1},\ldots ,s_{r}) \\
&&\times \,g(t_{p-r+1},\ldots ,t_{p+q-2r},s_{1},\ldots ,s_{r})d\mu
(s_{1})\ldots d\mu (s_{r}).
\end{eqnarray*}
It can then be also shown that the following {\sl multiplication formula} holds: if $f\in \EuFrak
H^{\odot p}$ and $g\in \EuFrak
H^{\odot q}$, then
\begin{eqnarray}\label{multiplication}
I_p(f) I_q(g) = \sum_{r=0}^{p \wedge q} r! {p \choose r}{ q \choose r} I_{p+q-2r} (f\widetilde{\otimes}_{r}g).
\end{eqnarray}
\smallskip

Let us now introduce some basic elements of the Malliavin calculus with respect
to the isonormal Gaussian process $X$. Let $\mathcal{S}$
be the set of all
cylindrical random variables of
the form
\begin{equation}
F=g\left( X(\phi _{1}),\ldots ,X(\phi _{n})\right) ,  \label{v3}
\end{equation}
where $n\geq 1$, $g:\mathbb{R}^{n}\rightarrow \mathbb{R}$ is an infinitely
differentiable function with compact support and $\phi _{i}\in \EuFrak H$.
The {\sl Malliavin derivative}  of $F$ with respect to $X$ is the element of
$L^2(\Omega ,\EuFrak H)$ defined as
\begin{equation*}
DF\;=\;\sum_{i=1}^{n}\frac{\partial g}{\partial x_{i}}\left( X(\phi
_{1}),\ldots ,X(\phi _{n})\right) \phi _{i}.
\end{equation*}
In particular, $DX(h)=h$ for every $h\in \EuFrak H$. By iteration, one can
define the $m$th derivative $D^{m}F$, which is an element of $L^2(\Omega ,\EuFrak H^{\odot m})$,
for every $m\geq 2$.
For $m\geq 1$ and $p\geq 1$, ${\mathbb{D}}^{m,p}$ denotes the closure of
$\mathcal{S}$ with respect to the norm $\Vert \cdot \Vert _{m,p}$, defined by
the relation
\begin{equation*}
\Vert F\Vert _{m,p}^{p}\;=\;E\left[ |F|^{p}\right] +\sum_{i=1}^{m}E\left(
\Vert D^{i}F\Vert _{\EuFrak H^{\otimes i}}^{p}\right) .
\end{equation*}
The Malliavin derivative $D$ verifies the following \textsl{chain rule}. If
$\varphi :\mathbb{R}^{n}\rightarrow \mathbb{R}$ is continuously
differentiable with bounded partial derivatives and if $F=(F_{1},\ldots
,F_{n})$ is a vector of elements of ${\mathbb{D}}^{1,2}$, then $\varphi
(F)\in {\mathbb{D}}^{1,2}$ and
\begin{equation*}
D\,\varphi (F)=\sum_{i=1}^{n}\frac{\partial \varphi }{\partial x_{i}}(F)DF_{i}.
\end{equation*}
Note also that a random variable $F$ as in (\ref{E}) is in ${\mathbb{D}}^{1,2}$ if and only if
$\sum_{q=1}^{\infty }q\|J_qF\|^2_{L^2(\Omega)}<\infty$
and, in this case, $E\left( \Vert DF\Vert _{\EuFrak H}^{2}\right)
=\sum_{q=1}^{\infty }q\|J_qF\|^2_{L^2(\Omega)}$. If $\EuFrak H=
\mathrm{L}^{2}(A,\mathcal{A},\mu )$ (with $\mu $ non-atomic), then the
derivative of a random variable $F$ as in (\ref{E}) can be identified with
the element of $L^2(A\times \Omega )$ given by
\begin{equation}
D_{x}F=\sum_{q=1}^{\infty }qI_{q-1}\left( f_{q}(\cdot ,x)\right) ,\quad x\in
A.  \label{dtf}
\end{equation}

We denote by $\delta $ the adjoint of the operator $D$, also called the
\textsl{divergence operator}. A random element $u\in L^2(\Omega ,\EuFrak H)$
belongs to the domain of $\delta $, noted $\mathrm{Dom}\delta $, if and
only if it verifies
$|E\langle DF,u\rangle _{\EuFrak H}|\leq c_{u}\,\Vert F\Vert _{L^2(\Omega)}$
for any $F\in \mathbb{D}^{1,2}$, where $c_{u}$ is a constant depending only
on $u$. If $u\in \mathrm{Dom}\delta $, then the random variable $\delta (u)$
is defined by the duality relationship (called \textsl{integration by parts
formula})
\begin{equation}
E(F\delta (u))=E\langle DF,u\rangle _{\EuFrak H},  \label{ipp}
\end{equation}
which holds for every $F\in {\mathbb{D}}^{1,2}$.
The divergence operator $\delta $ is also called the
\textsl{Skorohod integral} because in the case of the Brownian motion it coincides
with the anticipating stochastic integral introduced by Skorohod in \cite{Sk}.

The family $(T_t,\,t\geq 0)$ of operators is defined through the projection operators $J_q$
as \begin{equation}\label{OUsemigroup}
T_t=\sum_{q=0}^\infty e^{-qt}J_q,
\end{equation}
and is called the {\sl Ornstein-Uhlenbeck semigroup}. Assume that
the process $X'$, which stands for an independent copy of $X$, is such that $X$ and $X'$ are defined
on the product probability space $(\Omega\times\Omega',\mathscr{F}\otimes\mathscr{F}',P\times P')$.
Given a random variable $Z\in\mathbb{D}^{1,2}$, we can regard its Malliavin derivative
$DZ=DZ(X)$ as a measurable mapping from
$\R^\HH$ to $\R$, determined $P\circ X^{-1}$-almost surely. Then, for any $t\geq 0$, we have the
so-called {\sl Mehler's formula} (see e.g. \cite[Section 8.5,
Ch. I]{MallBook} or \cite[formula
(1.54)]{Nbook}):
\begin{equation}\label{mehler}
T_t(DZ)=E'\big(DZ(e^{-t}X+\sqrt{1-e^{-2t}}X')\big),
\end{equation}
where $E'$ denotes the mathematical expectation with respect to the probability $P'$.

The operator $L$ is defined as
$L=\sum_{q=0}^{\infty }-qJ_{q}$,
and it can be proven to be the infinitesimal generator of the Ornstein-Uhlenbeck
semigroup $(T_t)_{t\geq 0}$.
The domain of $L$ is
\begin{equation*}
\mathrm{Dom}L=\{F\in L^2(\Omega ):\sum_{q=1}^{\infty }q^{2}\left\|
J_{q}F\right\| _{L^2(\Omega )}^{2}<\infty \}=\mathbb{D}^{2,2}\text{.}
\end{equation*}
There is an important relation between the operators $D$, $\delta $ and $L$
(see e.g. \cite[Proposition 1.4.3]{Nbook}). A random variable $F$ belongs to
$\mathbb{D}^{2,2}$ if and only if $F\in \mathrm{Dom}\left( \delta D\right) $
(i.e. $F\in {\mathbb{D}}^{1,2}$ and $DF\in \mathrm{Dom}\delta $), and in
this case
\begin{equation}
\delta DF=-LF.  \label{k1}
\end{equation}

For any $F \in L^2(\Omega )$, we define $L^{-1}F =\sum_{q=0}^{\infty }-\frac{1}{q} J_{q}(F)$. The operator $L^{-1}$ is called the
\textsl{pseudo-inverse} of $L$. For any $F \in L^2(\Omega )$, we have that $L^{-1} F \in  \mathrm{Dom}L$,
and
\begin{equation}\label{Lmoins1}
LL^{-1} F = F - E(F).
\end{equation}

We end the preliminaries by noting that Shigekawa \cite{shigekawa}
has developed an alternative framework which avoids the inverse of
the Ornstein-Uhlenbeck operator $L$. This framework could provide an
alternative derivation of the integration by parts formula (2.30) in
\cite{stein-ptrf} which leads to Theorem \ref{MallSteinBound}.

\section{Poincaré-type inequalities and bounds on distances}\label{SS : TechLemmas}
\subsection{Poincaré inequalities} \label{SS : 1Poinc}

The following statement contains, among others, a general version
(\ref{poincare1}) of the infinite-dimensional Poincaré inequality (\ref{INFpoinc}).
\begin{proposition}\label{nicelemma}
Fix $p\geq 2$ and let $F\in\mathbb{D}^{1,p}$ be such that $E(F)=0$.
\begin{enumerate}
\item The following estimate holds:
\begin{equation}\label{DF}
E\left\Vert DL^{-1}F\right\Vert _{\mathfrak{H}}^{p}\leq E\left\Vert
DF\right\Vert _{
\mathfrak{H}
}^{p}\text{.}
\end{equation}
\item If in addition $F\in \mathbb{D}^{2,p}$, then
\begin{equation}\label{D^2F}
E\left\Vert D^{2}L^{-1}F\right\Vert _{op}^{p}\leq \frac{1}{2^p}\,E\left\Vert
D^{2}F\right\Vert _{op}^{p},
\end{equation}
where $\left\Vert
D^{2}F\right\Vert _{op}$ indicates the operator norm of the
random Hilbert-Schmidt operator
\begin{equation*}
\mathfrak{H}\rightarrow \mathfrak{H} : f\mapsto \left\langle f,D^{2}F\right\rangle _{\mathfrak{H}}\text{.}
\end{equation*}
(and similarly for $\| D^2L^{-1}F \|_{op}$).
\item If $p$ is an even integer, then
\begin{equation}\label{poincare1}
E\big[F^{p}\big]\leq (p-1)^{p/2}\, E\big[\|DF\|_\HH^{p}\big].
\end{equation}
\end{enumerate}
\end{proposition}
{\it Proof}. By virtue of standard arguments,
we may assume throughout the proof that $\EuFrak H=\mathrm{L}^{2}(A,\mathcal{A},\mu )$, where
$(A,\mathcal{A})$ is a measurable space and $\mu $ is a $\sigma $-finite and
non-atomic measure.
\begin{enumerate}
\item In what follows, we will write $X'$ to indicate an independent copy of $X$. Let $F \in L^2(\Omega )$
have the expansion \eqref{E}. Then, from \eqref{dtf},
\begin{equation*}
-D_{x}L^{-1}F=\sum_{q\geq 1}I_{q-1}\left( f_q\left( x,\cdot \right) \right) \text{,}.
\end{equation*}
By combining this  relation
with Mehler's formula (\ref{mehler}), one deduces that
\begin{eqnarray*}
-D_{x}L^{-1}F &=&\int_{0}^{\infty }e^{-t}T_{t}D_{x}F(X) dt
=\int_{0}^{\infty }e^{-t}\,E_{X'}\, D_{x}F\left( e^{-t}X+\sqrt{1-e^{-2t}}
X'\right) dt\\
&=&E_{Y}\,E_{X'}\, D_{x}F\left( e^{-Y}X+\sqrt{1-e^{-2Y}}X'\right)
\end{eqnarray*}
where $Y\sim\mathcal{E}(1)$ is an
independent exponential random variable of mean $1$,
and $\{T_t : t\geq 0\}$ is the Ornstein-Uhlenbeck semigroup (\ref{OUsemigroup}). Note that we regard every random
variable $D_{x}F$ as an application $\mathbb{R}^{\mathfrak{H}}\rightarrow
\mathbb{R}$ and that (for a generic random variable $G$) we write $E_{G}$ to indicate that
we take the expectation with respect to $G$. It follows that
\begin{eqnarray*}
E\left\Vert DL^{-1}F\right\Vert _{\mathfrak{H}}^{p}
&=& E_X\left\|
E_{Y}\,E_{X'}\, DF\left( e^{-Y}X+\sqrt{1-e^{-2Y}}X'\right)
 \right\|_{\mathfrak{H}}^{p}\\
&\leq & E_X\,
E_{Y}\,E_{X'}\left\| DF\left( e^{-Y}X+\sqrt{1-e^{-2Y}}X'\right)
 \right\|_{\mathfrak{H}}^{p}\\
&=& E_Y\,
E_{X}\,E_{X'}\left\| DF\left( e^{-Y}X+\sqrt{1-e^{-2Y}}X'\right)
 \right\|_{\mathfrak{H}}^{p}\\
&= & E_Y\,
E_{X}\left\| DF\left( X\right)\right\|_{\mathfrak{H}}^{p}=
E_{X}\left\| DF\left( X\right)\right\|_{\mathfrak{H}}^{p}=
E\left\| DF\right\|_{\mathfrak{H}}^{p}
\end{eqnarray*}
where we used the fact
that $e^{-t}X'+\sqrt{1-e^{-2t}}X\overset{law}{=}X$ for any $t\geq 0$.
\item From the relation
\begin{equation*}
-D_{xy}^{2}L^{-1}F=\sum_{q\geq 2}\left( q-1\right) I_{q-2}\left( f_q\left( x,y,\cdot
\right) \right) \text{}
\end{equation*}
one deduces analogously that
\begin{eqnarray*}
-D_{xy}^{2}L^{-1}F &=&\int_{0}^{\infty }e^{-2t}T_{t}D_{xy}^{2}Fdt \\
&=&\int_{0}^{\infty }e^{-2t}\,E_{X'}\, D_{xy}^{2}F\left( e^{-t}X+\sqrt{1-e^{-2t}
}X'\right)dt\\
&=&\frac12\,E_Y\,E_{X'}\, D_{xy}^{2}F\left( e^{-Y}X+\sqrt{1-e^{-2Y}}X'\right)
\end{eqnarray*}
where $Y\sim\mathcal{E}(2)$ is an
independent exponential random variable of mean $\frac12$. Thus
\begin{eqnarray*}
E\left\Vert D^{2}L^{-1}F\right\Vert _{op}^{p}
&=&\frac1{2^p}E_{X}\left\| E_{Y}\,E_{X'}\,D^{2}F\left( e^{-Y}X+\sqrt{1-e^{-2Y}}X'\right)\right\|^p_{op}\\
&\leq&\frac1{2^p}E_{X}\,E_{Y}\,E_{X'}\,\left\|D^{2}F\left( e^{-Y}X+\sqrt{1-e^{-2Y}}X'\right)\right\|^p_{op}\\
&=&\frac1{2^p}E_{Y}\,E_{X}\,E_{X'}\,\left\|D^{2}F\left( e^{-Y}X+\sqrt{1-e^{-2Y}}X'\right)\right\|^p_{op}\\
&=&\frac1{2^p}E_{Y}\,E_{X}\,\left\|D^{2}F\left( X\right)\right\|^p_{op}
=\frac1{2^p}E_{X}\,\left\|D^{2}F\left( X\right)\right\|^p_{op}
=\frac1{2^p}E\left\|D^{2}F\right\|^p_{op}.
\end{eqnarray*}
\item Writing $p=2k$, we have
\begin{eqnarray*}
E\big[F^{2k}\big]&=&E\big[LL^{-1}F\times F^{2k-1}\big]=-E\big[\delta DL^{-1}F\times F^{2k-1}\big]\\
&=&(2k-1)E\big[\langle DF,-DL^{-1}F\rangle F^{2k-2}\big]\\
&\leq&(2k-1)\left(E\big[\big|\langle DF,-DL^{-1}F\rangle\big|^k\big]\right)^{\frac1k}
\left(E\big[F^{2k}\big]\right)^{1-\frac{1}{k}}\quad\mbox{by H\"older's  inequality},
\end{eqnarray*}
from which we infer that
\begin{eqnarray*}
E\big[F^{2k}\big]&\leq&(2k-1)^k E\big[\big|\langle DF,-DL^{-1}F\rangle\big|^k\big] \leq (2k-1)^k E\big[\| DF\|_\HH^k \|DL^{-1}F\|_\HH^k\big]\\
&\leq &(2k-1)^k \sqrt{E\big[\| DF\|_\HH^{2k}\big]}\sqrt{E\big[ \|DL^{-1}F\|_\HH^{2k}\big]}\leq
(2k-1)^k E\big[\| DF\|_\HH^{2k}\big].
\end{eqnarray*}
\qed
\end{enumerate}

We also state the following technical result which will be needed in Section \ref{S : MAIN}. The proof is standard and omitted.

\begin{lemma}\label{L : DDD} Let $F$ and $G$ be two elements of $\mathbb{D}^{2,4}$. Then, the two random elements
$\langle D^2F,DG \rangle_{\mathfrak H}$ and $\langle DF,D^2G \rangle_{\mathfrak H}$ belong to $L^2(\Omega ,\EuFrak H)$.
Moreover, $\langle DF, DG \rangle_{\mathfrak H } \in \mathbb{D}^{1,2}$ and
\begin{equation}\label{chose}
D \langle DF, DG \rangle_{\mathfrak H } = \langle D^2F,DG \rangle_{\mathfrak H} + \langle DF,D^2G \rangle_{\mathfrak H}.
\end{equation}
\end{lemma}

\subsection{Bounds on the total variation and Wasserstein distances}\label{SS :distances}

Let $U,Z$ be two generic real-valued random variables. We recall that the {\sl total variation distance}
between the law of $U$ and the law of $Z$ is defined as
\begin{equation}\label{TV dst}
d_{TV}(U,Z) = \sup_A |P(U\in A) - P(Z\in A)|,
\end{equation}
where the supremum is taken over all Borel subsets $A$ of $\mathbb{R}$. For two random vectors $U$ and $Z$
with values in $\mathbb{R}^d$, $d\geq 1$, the {\sl Wasserstein distance}
between the law of $U$ and the law of $Z$ is
\begin{equation}\label{W dst}
d_{W}(U,Z) = \sup_{f : \|f\|_{Lip} \leq 1} |E[f(U)] - E[f(Z)]|,
\end{equation}
where $\|\cdot \|_{Lip}$ stands for the usual Lipschitz seminorm. We stress that the topologies
induced by $d_{TV}$ and $d_W$, on the class of all probability measures on $\mathbb{R}$, are strictly stronger
than the topology of weak convergence. The following statement has been proved in \cite[Theorem 3.1]{stein-ptrf} by means of Stein's method.

\begin{theorem}\label{MallSteinBound} Suppose that $Z\sim{\mathscr N}(0,1)$. Let $F\in\mathbb{D}^{1,2}$ and $E(F)=0$. Then,
\begin{equation}\label{Wbound}
d_{W}(F,Z) \leq E|1-\langle DF, -DL^{-1}F\rangle_{\mathfrak H}|\leq E[(1-\langle DF, -DL^{-1}F\rangle_{\mathfrak H})^2]^{1/2}.
\end{equation}
If moreover $F$ has an absolutely continuous distribution, then
\begin{equation}\label{TVbound}
d_{TV}(F,Z) \leq 2E|1-\langle DF, -DL^{-1}F\rangle_{\mathfrak H}|\leq 2 E[(1-\langle DF, -DL^{-1}F\rangle_{\mathfrak H})^2]^{1/2}.
\end{equation}
\end{theorem}

\section{Proof of Theorem \ref{thm1} and contraction inequalities}\label{S : MAIN}
\subsection{Proof of Theorem \ref{thm1}} \label{SS :ProofT1}
We can assume, without loss of generality, that $\mu=0$ and $\sigma^2=1$. Set
$W=\left\langle DF,-DL^{-1}F\right\rangle _{\mathfrak{H}}$. First, note that $W$ has mean 1, as
$$
E(W)=E[\langle DF,-DL^{-1}F\rangle_{\mathfrak{H}}]= - E[F \times \delta D L^{-1} F ] = E[F \times L L^{-1} F]
= E [F^2] =1.
$$
By Theorem \ref{MallSteinBound} it follows
that we only need to bound $\sqrt{{\rm Var}\left( W \right)}$.
By (\ref{INFpoinc}), we have
${\rm Var}\left( W \right) \leq E \left\Vert D W \right\Vert _{\mathfrak{H}}^{2}$.
So, our problem is now to evaluate $\left\Vert D W \right\Vert _{\mathfrak{H}}^{2}$. By using Lemma \ref{L : DDD} in the special case $G = -L^{-1}F$, we deduce that
\begin{eqnarray*}
\left\Vert DW \right\Vert _{\mathfrak{H}}^{2} &=&
\left\| \langle D^2F,-DL^{-1}F\rangle_\mathfrak{H} + \langle DF,-D^2L^{-1}F\rangle_\mathfrak{H}
\right\|_\mathfrak{H}^2 \\
&\leq &
2\left\|\langle D^2F,-DL^{-1}F\rangle_\mathfrak{H} \right\|^2_\mathfrak{H}
+2\left\|\langle DF,-D^2L^{-1}F\rangle_\mathfrak{H}\right\|^2_\mathfrak{H}.
\end{eqnarray*}
We evaluate the last two terms separately. We have
\begin{equation*}
\left\|\langle D^2F,-DL^{-1}F\rangle_\mathfrak{H} \right\|^2_\mathfrak{H}
\leq \left\Vert D^{2}F\right\Vert _{op}^{2}\left\Vert
DL^{-1}F\right\Vert _{\mathfrak{H}}^{2}
\end{equation*}
and
\begin{equation*}
\left\|\langle DF,-D^2L^{-1}F\rangle_\mathfrak{H}\right\|^2_\mathfrak{H}
\leq \left\Vert DF\right\Vert _{\mathfrak{H}}^{2}\left\Vert
D^{2}L^{-1}F\right\Vert _{op}^{2}\text{.}
\end{equation*}
It follows that
\begin{eqnarray*}
E \left\Vert DW \right\Vert _{\mathfrak{H}}^{2}
&\leq &2\,E\left[ \left\Vert DL^{-1}F\right\Vert _{\mathfrak{H}}^{2}\left\Vert D^{2}F\right\Vert _{op}^{2}
+\left\Vert DF\right\Vert_{\mathfrak{H}}^{2}\left\Vert D^{2}L^{-1}F\right\Vert _{op}^{2}\right] \\
&\leq &2\, \left( E\left\Vert DL^{-1}F\right\Vert _{\mathfrak{H}}^{4}\times E\left\Vert D^{2}F\right\Vert _{op}^{4}\right) ^{1/2}+2\,\left(
E\left\Vert DF\right\Vert _{\mathfrak{H}}^{4}\times E\left\Vert
D^{2}L^{-1}F\right\Vert _{op}^{4}\right) ^{1/2} \text{.}
\end{eqnarray*}
The desired conclusion follows by using, respectively, (\ref{DF}) and (\ref{D^2F}) with $p=4$.
\qed

\vskip.5cm

\subsection{Random contraction inequalities}\label{SS : CONTR}

When the quantity $E\left\Vert D^{2}F\right\Vert _{op}^{4}$ appearing in (\ref{Poincw})-(\ref{Poinc}) is analytically
too hard to assess, one can resort to the following inequality, which  we name
\textsl{random contraction inequality}:

\begin{proposition}{\bf (Random contraction inequality)}. \label{T : RandomContr}
Let $F\in \mathbb{D}^{2,4}$. Then
\begin{equation}\label{INOP}
\left\Vert D^{2}F \right\Vert
_{op}^{4}\leq \left\Vert D^{2}F \otimes
_{1}D^{2}F \right\Vert _{\mathfrak{H}^{\otimes 2}}^{2}
\text{,}
\end{equation}
where $D^{2}F\otimes _{1}D^{2}F$ is the random element of $\mathfrak{H}^{\odot 2}$ obtained as the
contraction of the symmetric random tensor $D^{2}F$, see (\ref{v2}).
\end{proposition}
{\it Proof}.
We can associate with the symmetric random elements $D^{2}F\in \mathfrak{H}^{\odot 2}$ the
random Hilbert-Schmidt operator
$\, f\mapsto \left\langle f,D^{2}F\right\rangle _{\mathfrak{H}^{\otimes 2}}\text{.}$
Denote by $\{ \gamma _{j}\}_{j\geq 1}$ the sequence of its (random) eigenvalues.
One has that
\begin{equation*}
\left\Vert D^{2}F \right\Vert
_{op}^{4}=\max_{j\geq 1}\left\vert \gamma _{j} \right\vert
^{4}\leq \sum_{j\geq 1}\left\vert \gamma _{j}
\right\vert ^{4}=\left\Vert D^{2}F \otimes
_{1}D^{2}F \right\Vert _{\mathfrak{H}^{\otimes 2}}^{2},
\end{equation*}
and the conclusion follows.
\qed

\medskip
The following result is an immediate corollary of Theorem \ref{thm1} and Proposition \ref{T : RandomContr}.

\begin{cor}
Let $F\in \mathbb{D}^{2,4}$ with $E\left( F\right) =\mu $ and ${\rm Var}\left( F\right) =\sigma ^{2}$.
Assume that $Z\sim\mathscr{N}( \mu ,\sigma ^{2}) $.
Then
\begin{equation}
d_{W}\left( F,Z \right) \leq \frac{\sqrt{10}}{2 \sigma}
E\left[ \Vert D^{2}F\otimes _{1}D^{2}F\Vert _{\mathfrak{H}^{\otimes 2}}^{2}\right] ^{\frac14}
\times E\left[ \left\Vert DF\right\Vert _{
\mathfrak{H}
}^{4}\right] ^{\frac14}.  \label{ContrW}
\end{equation}
If, in addition,
the law of $F$ is absolutely continuous with respect to the Lebesgue measure, then
\begin{equation}
d_{TV}\left( F,Z \right) \leq \frac{\sqrt{10}}{\sigma ^{2}}
E\left[ \Vert D^{2}F\otimes _{1}D^{2}F\Vert _{\mathfrak{H}^{\otimes 2}}^{2}\right] ^{\frac14}
\times E\left[ \left\Vert DF\right\Vert _{
\mathfrak{H}
}^{4}\right] ^{\frac14}.   \label{Contr}
\end{equation}
\end{cor}

\medskip
\begin{remark}\label{rkgio}
{\rm
When used in the context of central limit theorems, inequality (\ref{Contr}) does not give, in general, optimal
rates. For instance, if $F_{k}=I_{2}\left( f_{k}\right) $ is a sequence of
double integrals such that $E\left( F_{k}^{2}\right) \rightarrow 1$ and
$F_{k}\overset{law}{\longrightarrow }Z\sim \mathscr{N}\left( 0,1\right)$ as $k\to\infty$, then (\ref{Contr})
implies that
\begin{equation*}
d_{TV}\left( F_{k},Z\right) \leq
cst\times\left\Vert f_{k}\otimes _{1}f_{k}\right\Vert _{\mathfrak{H}^{\otimes
2}}^{1/2}\rightarrow 0\text{,}
\end{equation*}
and the rate $\left\Vert f_{k}\otimes _{1}f_{k}\right\Vert _{\mathfrak{H}^{\otimes 2}}^{1/2}$ is suboptimal (by a power of $1/2$), see Proposition 3.2 in \cite{stein-ptrf}.
}
\end{remark}

\section{Characterization of CLTs on a fixed Wiener chaos}\label{S : CharCLTWC}

The following statement collects results proved in \cite{NuPec} (for the equivalences between (i), (ii)
and (iii))
and \cite{NO} (for the equivalence with (iv)).

\begin{theorem}\label{T : CLT on WienerC} Fix $q\geq 2$, and let $F_{k}=I_{q}\left( f_{k}\right)$, $k\geq 1$, be a sequence of
multiple Wiener-It\^o integrals such that $E\left( F_{k}^{2}\right) \rightarrow 1$. As $k\to\infty$, the following
four conditions are equivalent:
\begin{description}
\item[\bf (i)] $F_k \stackrel {law}{\longrightarrow} Z\sim\mathscr{N}(0,1)$;

\item[\bf (ii)] $E(F_k^4)\longrightarrow E(Z^4)= 3$;

\item[\bf (iii)] $\|f_k\otimes_r f_k\|_{\HH^{\otimes(2q-2r)}}\longrightarrow 0$ for all $r=1,\ldots,q-1$;

\item[\bf (iv)] $\|DF\|^2_{\mathfrak H} \overset{L^2(\Omega)}{\longrightarrow} q $.
\end{description}
\end{theorem}

\bigskip

See Section 9 in \cite{PecTaqSurvey} for a discussion of the combinatorial aspects of the implication (ii) $\rightarrow$ (i)
in the statement of Theorem \ref{T : CLT on WienerC}. The next theorem, which is a consequence of the main results of this
paper, provides two new necessary and sufficient conditions for CLTs on a fixed Wiener chaos.

\begin{theorem}
Fix $q\geq 2$, and let $F_{k}=I_{q}\left( f_{k}\right) $ be a sequence of
multiple Wiener-It\^o integrals such that $E\left( F_{k}^{2}\right) \rightarrow 1$.
Then, the following three conditions are equivalent as $k\to\infty$:

\begin{description}

\item[(i)] $F_{k}\overset{law}{\longrightarrow }Z\sim\mathscr{N}(0,1)$;

\item[(ii)] $\left\Vert D^{2}F_{k}\otimes _{1}D^{2}F_{k}\right\Vert _{\mathfrak{H}^{\otimes 2}} \overset{L^2(\Omega)}{\longrightarrow } 0$;

\item[(iii)] $\left\Vert D^{2}F_{k}\right\Vert _{op} \overset{L^4(\Omega)}{\longrightarrow } 0$.
\end{description}
\end{theorem}
{\it Proof}.
Since $E\|DF_k\|^2_{\mathfrak{H}} = q E(F^2_k)\rightarrow q$, and since the
random variables $\|DF_k\|^2_{\mathfrak{H}}$ live inside a finite sum of Wiener chaoses
(where all the $L^p(\Omega)$ norms are equivalent), we deduce that the sequence $E\|DF_k\|^4_{\mathfrak{H}}$, $k\geq 1$,
is bounded. In view of (\ref{Poincw}) and (\ref{INOP}), it is therefore enough to prove the implication (i) $\rightarrow$ (ii). Without loss of generality, we can assume that $\HH=L^2(A,\mathscr{A},\mu)$ where $(A,\mathscr{A})$ is
a measurable space and $\mu$ is a $\sigma$-finite measure with no atoms.
Now observe that
$$
D^2_{a,b}F_k=q(q-1)I_{q-2}\big(f_k(\cdot,a,b)\big),\quad a,b\in A.
$$
Hence, using the multiplication formula \eqref{multiplication},
\begin{eqnarray*}
&&D^2F_k\otimes_1 D^2F_k(a,b)\\
&=&q^2(q-1)^2 \int_{A}I_{q-2}\big(f_k(\cdot,a,u)\big)
I_{q-2}\big(f_k(\cdot,b,u)\big)\mu(du)\\
&=&q^2(q-1)^2 \sum_{r=0}^{q-2}r!\binom{q-2}{r}^2
I_{2q-4-2r}\left(
\int_{A}f_k(\cdot,a,u)\widetilde{\otimes }_r
f_k(\cdot,b,u)\mu(du)\right)\\
&=&q^2(q-1)^2 \sum_{r=0}^{q-2}r!\binom{q-2}{r}^2
I_{2q-4-2r}\left(
f_k(\cdot,a)\widetilde{\otimes }_{r+1}
f_k(\cdot,b)\right)\\
&=&q^2(q-1)^2 \sum_{r=1}^{q-1}(r-1)!\binom{q-2}{r-1}^2
I_{2q-2-2r}\left(
f_k(\cdot,a)\widetilde{\otimes }_{r}
f_k(\cdot,b)\right).
\end{eqnarray*}
Using the orthogonality and isometry properties of the integrals $I_q$, we get
$$
E\left\Vert D^{2}F_{k}\otimes _{1}D^{2}F_{k}\right\Vert _{\mathfrak{H}^{\otimes 2}}^2
\leq q^4(q-1)^4 \sum_{r=1}^{q-1}(r-1)!^2\binom{q-2}{r-1}^4
(2q-2-2r)!
\|
f_k\otimes_{r}
f_k\|^2_{\HH^{\otimes(2q-2r)}}.
$$
The desired conclusion now follows since, according to
Theorem \ref{T : CLT on WienerC}, if (i) is verified then, necessarily,
$\|f_k\otimes_{r}
f_k\|_{\HH^{\otimes(2q-2r)}}\rightarrow 0$ for every $r=1,...,q-1$. \qed

\bigskip

\section{CLTs for linear functionals of Gaussian subordinated fields}\label{S : GaussSub}

We now provide an explicit application of the inequality (\ref{ContrW}).
Let $B$ denote a centered Gaussian process with stationary increments and such that $\int_\R|\rho(x)|dx<\infty$,
where $\rho(u-v):=E\big[(B_{u+1}-B_u)(B_{v+1}-B_v)\big]$. Also, in order to avoid trivialities, assume that
$\rho$ is not identically zero.

The Gaussian space generated by $B$ can be identified with an isonormal Gaussian process of the type
$X=\{X(h),\,h\in\HH\}$, for $\HH$ defined as follows: (i) denote by $\mathscr{E}$ the set of all
step functions on $\R$, (ii) define $\HH$ as the Hilbert space obtained by closing
$\mathscr{E}$ with respect to the inner product $\langle {\bf 1}_{[s,t]},{\bf 1}_{[u,v]}\rangle_\HH=
{\rm Cov}(B_t-B_s,B_v-B_u)$.
In particular, with such a notation, one has that $B_t-B_s=X({\bf 1}_{[s,t]})$.

Let $f:\R\to\R$ be a
real function of class $\mathscr{C}^2$, and $Z \sim\mathscr{N}(0,1)$.
We assume that $f$ is not constant, that $E|f(Z )|<\infty$ and that $E|f''(Z )|^4<\infty$.
As a consequence of the generalized
Poincar\'e inequality (\ref{poincare1}), we see that we also
automatically have $E|f'(Z)|^4<\infty$ and $E|f(Z)|^4<\infty$.

Fix $a<b$ in $\R$ and, for any $T>0$,
consider
$$
F_T = \frac{1}{\sqrt{T}}\int_{aT}^{bT} \big(f(B_{u+1}-B_u)-E[f(Z)]\big)du.
$$

\begin{thm}\label{tvtv}
As $T\to\infty$,
\begin{equation}\label{tvtvtv}
d_{W}\left(\frac{F_T}{\sqrt{{\rm Var}F_T}},Z\right)=O(T^{-1/4}).
\end{equation}
\end{thm}
\begin{rem}
{\rm
We believe that the rate in (\ref{tvtvtv}) is not optimal (it should be $O(T^{-1/2})$ instead), see also
Remark \ref{rkgio}.
}
\end{rem}
{\it Proof of Theorem \ref{tvtv}}. We have
$$
DF_T = \frac{1}{\sqrt{T}}\int_{aT}^{bT} f'(B_{u+1}-B_u){\bf 1}_{[u,u+1]}du
$$
and
$$
D^2F_T = \frac{1}{\sqrt{T}}\int_{aT}^{bT} f''(B_{u+1}-B_u){\bf 1}_{[u,u+1]}^{\otimes 2}du.
$$
Hence
$$
\|DF_T\|_{\mathfrak H}^2=\frac1T\int_{[aT,bT]^2} f'(B_{u+1}-B_u)f'(B_{v+1}-B_v)\,\rho(u-v)dudv
$$
so that
\begin{eqnarray*}
\|DF_T\|_{\mathfrak H}^4&=&\frac1{T^2}\int_{[aT,bT]^4} f'(B_{u+1}-B_u)f'(B_{v+1}-B_v)
f'(B_{w+1}-B_w)\\
&&\hskip3cm \times f'(B_{z+1}-B_z)\rho(w-z)\rho(u-v)
dudvdwdz.
\end{eqnarray*}
By applying Cauchy-Schwarz inequality twice, and by using the fact that $B_{u+1}-B_u \overset{law}{=} Z$, we get
$$
\left|E\big(f'(B_{u+1}-B_u)f'(B_{v+1}-B_v)
f'(B_{w+1}-B_w)f'(B_{z+1}-B_z)\big)\right|
\leq
E|f'(Z)|^4
$$
so that
\begin{eqnarray}
E\|DF_T\|_{\mathfrak H}^4 &\leq& E|f'(Z)|^4\left(\frac1T\int_{[aT,bT]^2} |\rho(u-v)|dudv\right)^2\notag\\
&\leq& E|f'(Z)|^4\left(\frac1T\int_{aT}^{bT} du \int_{\R} |\rho(x)|dx\right)^2 =O(1).\label{jenaibesoin}
\end{eqnarray}
On the other hand, we have
$$
D^2F_T\otimes_1 D^2F_T = \frac1T\int_{[aT,bT]^2}f''(B_{u+1}-B_u)f''(B_{v+1}-B_v)\rho(u-v){\bf 1}_{[u,u+1]}\otimes
{\bf 1}_{[v,v+1]} dudv.
$$
Hence
\begin{eqnarray*}
&&E\|D^2F_T\otimes_1 D^2F_T\|_{\mathfrak{H}^{\otimes 2}}^2\\
&=& \frac1{T^2}\int_{[aT,bT]^4}E\left(f''(B_{u+1}-B_u)f''(B_{v+1}-B_v)f''(B_{w+1}-B_w)
f''(B_{z+1}-B_z)\right)\\
&&\hskip1cm \times\rho(u-v)\rho(w-z)\rho(u-w)\rho(z-v) dudvdwdz\\
&\leq& E|f''(Z)|^4\frac1{T^2}\int_{[aT,bT]^4} |\rho(u-v)|
|\rho(w-z)||\rho(u-w)||\rho(z-v)|
dudvdwdz\\
&\leq& E|f''(Z)|^4\,\frac{b-a}T \int_{\R^3} |\rho(x)||\rho(y)||\rho(t)||\rho(x-y-t)|dxdydt =O(T^{-1}).
\end{eqnarray*}
By combining all these facts and (\ref{ContrW}),
the desired conclusion follows. \qed

\medskip
Theorem \ref{tvtv} does not guarantee that $\lim_{T\to\infty}{\rm Var}F_T$ exists. The following
proposition shows that the limit does indeed exist, at least when $f$ is symmetric.

\begin{proposition} \label{tvtvsymm}
Suppose that $f:\R\to\R$ is a {\sl symmetric}
real function of class $\mathscr{C}^2$. Then $\sigma^2 := \lim_{T\to\infty}{\rm Var}F_T$ exists in $(0,\infty)$. Moreover, as $T
\rightarrow \infty$,
\begin{equation}\label{acompareravecmarcusetrosen}
F_T \overset{law}{\longrightarrow} Z\sim\mathscr{N}(0,\sigma^2).
\end{equation}
\end{proposition}
{\it Proof of Proposition \ref{tvtvsymm}}.
We expand $f$ in terms of Hermite polynomials. Since $f$ is symmetric, we can write
$$
f(x)=E[f(Z)]+\sum_{q=1}^{\infty} c_{2q}H_{2q}(x),\quad x\in\R,
$$
where the real numbers $c_{2q}$ are given by $(2q)!c_{2q}=E[f(Z)H_{2q}(Z)]$.
Thus
\begin{eqnarray*}
{\rm Var}F_T&=&\frac1T\int_{[aT,bT]^2} {\rm Cov}\big(f(B_{u+1}-B_u),f(B_{v+1}-B_v)\big)dudv\\
&=&\sum_{q=1}^\infty c_{2q}^2 (2q)! \,\frac1T\int_{[aT,bT]^2}\rho^{2q}(v-u) dudv\\
&=&\sum_{q=1}^\infty c_{2q}^2 (2q)! \,\frac1T\int_{aT}^{bT}du\int_{aT-u}^{bT-u}dx\rho^{2q}(x)\\
&=&\sum_{q=1}^\infty c_{2q}^2 (2q)! \,\int_a^bdu\int_{-T(u-a)}^{T(b-u)}dx\rho^{2q}(x)\\
&\underset{T\nearrow\infty}{\longrightarrow}&(b-a)\sum_{q=1}^\infty c_{2q}^2 (2q)! \,\int_{\R}
\rho^{2q}(x)dx=:\sigma^2,\quad\mbox{by monotone convergence}.
\end{eqnarray*}

Since $f$ is not constant, there exists $q\geq 1$ such that $c_{2q}\neq 0$ so that $\sigma^2>0$ (recall that
we assumed $\rho\not\equiv 0$).
Moreover, we also have
$$
{\rm Var}F_T\leq E\big[\|DF_T\|_{\mathfrak H}^2\big]\leq \sqrt{E\big[\|DF_T\|_{\mathfrak H}^4\big]}=O(1), \quad\mbox{see (\ref{jenaibesoin})},
$$
so that $\sigma^2<\infty$. The assertion now follows from
 Theorem \ref{tvtv}.
\qed

\vskip.5cm

When $B$ is a fractional Brownian motion with Hurst index $H<1/2$, Theorem \ref{tvtv} applies because,
in this case, it is easily checked that $\int_\R|\rho(x)|dx<\infty$. On the other hand, using the
scaling property of $B$, observe that $F_{1/h}\overset{law}{=}\frac{1}{\sqrt{h}}
\int_a^b \left[f\left(\frac{B_{x+h}-B_x}{h^H}\right) -
E\big(f(Z)\big)\right] dx$ for all fixed $h>0$. Hence, since $E|B_t-B_s|^2=\sigma^2(|t-s|)$ with
$\sigma^2(r)=r^{2H}$ a concave function, the general Theorem 1.1 in \cite{MR} also applies, and this gives
another proof of (\ref{acompareravecmarcusetrosen}). We believe however that,
even in this particular case,
our proof is simpler (since not based on the rather technical method of moments). Moreover, note that \cite{MR}
is not concerned with bounds on distance between the laws of $F_{1/h}/\sqrt{{\rm Var}F_{1/h}}$ and
$Z\sim\mathscr{N}(0,1)$.

\section{A multidimensional extension}\label{S : Generaliz}
Let $V,Y$ be two random vectors with values in $\mathbb{R}^d$, $d\geq 2$. Recall that the Wasserstein distance
between the laws of $V$ and $Y$ is defined in (\ref{W dst}).
 The following statement, whose proof is based on the
 results obtained in \cite{npr-ihp}, provides a multidimensional version of (\ref{Poincw}).

\begin{thm}
Fix $d\geq 2$, and let $C=\{C(i,j):\,i,j=1,\ldots,d\}$ be a $d\times d$ positive definite matrix.
Suppose that $F=(F_1,\ldots,F_d)$ is a $\R^d$-valued random vector
such that $E[F_i]=0$ and $F_i\in\mathbb{D}^{2,4}$ for every $i=1,\ldots,d$.
Assume moreover that $F$ has covariance matrix $C$. Then
$$
d_W\big(F,\mathscr{N}_d(0,C)\big)\leq \frac{3\sqrt{2}}{2}\|C^{-1}\|_{op}\|C\|_{op}^{1/2}
\sum_{i=1}^d \big( E\|D^2 F_i\|^4_{op}\big)^{1/4}\times\sum_{j=1}^d\big(E\|DF_j\|^4_{\HH}\big)^{1/4},
$$
where $\mathscr{N}_d(0,C)$ indicates a $d$-dimensional centered Gaussian vector, with covariance matrix equal to $C$.
\end{thm}
{\it Proof}.
In \cite[Theorem 3.5]{npr-ihp} it is shown that
$$
d_W\big(F,\mathscr{N}_d(0,C)\big)
\leq \|C^{-1}\|_{op}\|C\|_{op}^{1/2}
\sqrt{
\sum_{i,j=1}^d E\big[(C(i,j)-\langle DF_i,-DL^{-1}F_j\rangle_\HH)^2\big]
} .
$$
Since, using successively (\ref{ipp}), (\ref{k1}) and (\ref{Lmoins1}), we have
$$E\big[\langle DF_i,-DL^{-1}F_j\rangle_\HH\big]
=-E\big[F_i\times \delta DL^{-1}F_j]=E\big[F_i\times LL^{-1}F_j]=E[F_iF_j]=C(i,j),$$
we deduce,
applying successively (\ref{INFpoinc}), (\ref{chose}), Cauchy-Schwarz inequality and Proposition \ref{nicelemma},
\begin{eqnarray*}
&&d_W\big(F,\mathscr{N}_d(0,C)\big)\\
&\leq&
\|C^{-1}\|_{op}\|C\|_{op}^{1/2}
\sum_{i,j=1}^d \sqrt{{\rm Var}\big[\langle DF_i,-DL^{-1}F_j\rangle_\HH\big]}\\
&\leq&
\|C^{-1}\|_{op}\|C\|_{op}^{1/2}
\sum_{i,j=1}^d \sqrt{E\big[
\|D\langle DF_i,-DL^{-1}F_j\rangle_\HH\|^2_\HH
\big]}\\
&\leq&
\sqrt{2}\|C^{-1}\|_{op}\|C\|_{op}^{1/2}
\sum_{i,j=1}^d \left(\sqrt{E\big[\|\langle D^2F_i,-DL^{-1}F_j\rangle_\HH\|^2_\HH\big]}
+\sqrt{E\big[\|\langle DF_i,-D^2L^{-1}F_j\rangle_\HH\|^2_\HH\big]}\right)\\
&\leq&
\sqrt{2}\|C^{-1}\|_{op}\|C\|_{op}^{1/2}
\sum_{i,j=1}^d
\left[
\left(E\big[\|D^2F_i\|^4_{op}\big]\right)^{1/4}
\left(E\big[\|DL^{-1}F_j\|^4_{\HH}\big]\right)^{1/4}
\right.\\
&&\left.\hskip7cm
+
\left(E\big[\|DF_i\|^4_{\HH}\big]\right)^{1/4}
\left(E\big[\|D^2L^{-1}F_j\|^4_{op}\big]\right)^{1/4}
\right]\\
&\leq&
\sqrt{2}\|C^{-1}\|_{op}\|C\|_{op}^{1/2}
\sum_{i,j=1}^d
\left[
\left(E\big[\|D^2F_i\|^4_{op}\big]\right)^{1/4}
\left(E\big[\|DF_j\|^4_{\HH}\big]\right)^{1/4}
\right.\\
&&\left.\hskip7cm
+
\frac12\left(E\big[\|DF_i\|^4_{\HH}\big]\right)^{1/4}
\left(E\big[\|D^2F_j\|^4_{op}\big]\right)^{1/4}
\right]\\
&=&
\frac{3\sqrt{2}}{2}\|C^{-1}\|_{op}\|C\|_{op}^{1/2}
\sum_{i=1}^d
\left(E\big[\|D^2F_i\|^4_{op}\big]\right)^{1/4}
\times\sum_{j=1}^d\left(E\big[\|DF_j\|^4_{\HH}\big]\right)^{1/4}.
\end{eqnarray*}
\qed

\vskip1cm

\noindent {\bf Acknowledgments.} We would like to thank Professor Paul Malliavin
for very stimulating discussions and an anonymous referee for pointing us towards the reference \cite{shigekawa}.


\begin{thebibliography}{99}

\bibitem{BBCG} D. Bakry, F. Barthe, P. Cattiaux and A. Guillin (2008). A simple proof of the Poincar\'{e} inequality for
a large class of probability measures including the log-concave case. {\it Elect. Comm. in Probab.} {\bf 13}, 60--66 (Electronic)

\bibitem{Bobkov} S.G. Bobkov (1999). Isoperimetric and analytic inequalities for log-concave probability measures.
{\it Ann. Probab.} {\bf 27}(4), 1903--1921.

\bibitem{CacPaUt} T. Cacoullos, V. Papahanasiou and S.A. Utev (1994). Variational inequalities with examples and an
 application to the central limit theorem. {\it Ann. Probab.} {\bf 22}(3), 1607--1618.

\bibitem{Chatterjee_ptrf} S. Chatterjee (2009). Fluctuation of eigenvalues and second order Poincaré inequalities.
 \textit{Probab. Theory Related Fields} {\bf 143}, 1--40.

\bibitem{Chen1985} L.H.Y. Chen (1985). Poincar\'{e}-type inequalities via stochastic integrals.
 {\it Z. Wahrsch. verw. Gebiete} {\bf 69}, 251--277.

\bibitem{Chen1987} L.H.Y. Chen (1987). Characterization of probability distributions by Poincar\'{e}-type inequalities.
{\it Ann. Inst. H. Poincar\'e Probab. Statist.} {\bf 23}(1), 91--110.

\bibitem{Chen1988} L.H.Y. Chen (1988). The central limit theorem and Poincar\'{e}-type inequalities. \textit{Ann. Probab.}
{\bf 16}(1), 300--304.

\bibitem{Chen_Shao_sur} L. Chen and Q.-M. Shao (2005).
Stein's method for normal approximation. In: \textit{An
introduction to Stein's method}, 1--59. Lect. Notes Ser. Inst.
Math. Sci. Natl. Univ. Singap. \textbf{4}, Singapore Univ. Press,
Singapore, 2005.

\bibitem{Chernoff} H. Chernoff (1981). A note on an inequality involving the normal distribution.
{\it Ann. Probab.} {\bf 9}(3), 533--535.

\bibitem{HouPA} C. Houdré and V. Pérez-Abreu (1995). Covariance
identities and inequalities for func- tionals on Wiener and
Poisson spaces. \textit{Ann. Probab.} \textbf{23}, 400--419.

\bibitem{It}
\rm K. It\^{o} (1951).
\rm Multiple Wiener integral.
{\it J. Math. Soc. Japan} \textbf{3}, 157--169.

\bibitem{MallBook} P. Malliavin (1997). \textit{Stochastic Analysis}.
Springer-Verlag, Berlin, Heidelberg, New York.

\bibitem{MR}
\rm M. Marcus and J. Rosen (2008).
\rm CLT for $L^p$ moduli of continuity of Gaussian processes.
{\it Stoch. Proc. Appl.} {\bf 118}, 1107--1135.


\bibitem{Nash} J. Nash (1956). Continuity of solutions of parabolic and elliptic equations. {\it Amer. J. Math.} {\bf 80},
 931--954.

\bibitem{NouPecNcltAOP}
\rm I. Nourdin and G. Peccati (2008).
\rm Non-central convergence of multiple integrals.
{\it Ann. Probab.}, to appear.

\bibitem{stein-ptrf}
\rm I. Nourdin and G. Peccati (2008).
\rm Stein's method on Wiener chaos.
{\it Probab. Theory Rel. Fields}, to appear.

\bibitem{NouPecexact} I. Nourdin and G. Peccati (2008).
Stein's method and exact Berry-Ess\'{e}n bounds for functionals of Gaussian fields. Preprint.

\bibitem{npr-ihp}
\rm I. Nourdin, G. Peccati and A. R\'eveillac (2008).
\rm Multivariate normal approximation using Stein's method and Malliavin calculus.
{\it Ann. Inst. H. Poincar\'e Probab. Statist.}, to appear.

\bibitem{Nbook}
\rm D. Nualart (2006).
\it The Malliavin calculus and related topics.
\rm Springer-Verlag, Berlin, 2nd edition.

\bibitem{NO}
\rm D. Nualart and S. Ortiz-Latorre (2008).
\rm Central limit theorem for multiple stochastic integrals and Malliavin calculus.
{\it Stoch. Proc. Appl.} {\bf 118} (4), 614--628.

\bibitem{NuPec}
\rm D. Nualart and G. Peccati (2005).
\rm Central limit theorems for sequence of multiple stochastic integrals.
{\it Ann. Probab.} {\bf 33} (1), 177--193.

\bibitem{PecTaqSurvey}
G. Peccati and M.S. Taqqu (2008). Moments, cumulants and diagram
formulae for non-linear functionals of random measures. Preprint available at \texttt{http://fr.arxiv.org/abs/0811.1726}.

\bibitem{PT}
G. Peccati and C.A. Tudor (2005). Gaussian limits for
vector-valued multiple stochastic integrals. In: \it S\'eminaire
de Probabilit\'es XXXVIII \rm, 247--262. Lecture Notes in Math.
{\bf 1857}, Springer-Verlag, Berlin.

\bibitem{Reinert_sur} G. Reinert (2005). Three general approaches to Stein's
method. In: \textit{An introduction to Stein's method}, 183--221.
Lect. Notes Ser. Inst. Math. Sci. Natl. Univ. Singap. \textbf{4},
Singapore Univ. Press, Singapore.

\bibitem{shigekawa}
\rm I. Shigekawa (1986).
\rm De Rham-Hodge-Kodaira's decomposition on an abstract Wiener space.
{\it J. Math. Kyoto Univ.} \textbf{26} (2), 191--202.



\bibitem{Sk}
\rm A. V. Skorohod (1975).
\rm On a generalization of a stochastic integral.
{\it Theory Probab. Appl.} \textbf{20}, 219--233.

\end{thebibliography}
\end{document}